\documentclass{proc-l}
\issueinfo{129}{10}{October}{2001}
\pagespan{3057}{3067}
\dateposted{April 16, 2001}
\PII{S 0002-9939(01)05928-7}
\copyrightinfo{2001}{American Mathematical Society}
\usepackage{amssymb}

\theoremstyle{plain}
\newtheorem{theorem}{Theorem}[section]
\newtheorem{corollary}[theorem]{Corollary}

\newtheorem{proposition}[theorem]{Proposition}

\theoremstyle{definition}

\theoremstyle{remark}
\newtheorem{remark}[theorem]{Remark}

\newcommand{\R}{\mathbb{R}}
\newcommand{\Rn}{\mathbb{R}^n}
\newcommand{\Z}{\mathbb{Z}}
\newcommand{\Zn}{\mathbb{Z}^n}
\newcommand{\Q}{\mathbb{Q}_A}
\newcommand{\QZ}{\mathbb{Q}_A \rtimes_{\vartheta} \mathbb{Z}}
\newcommand{\Ind}{\textstyle{\operatornamewithlimits{Ind}}}
\newcommand{\abs}[1]{\lvert #1\rvert}
\newcommand{\ip}[2]{\langle #1, #2\rangle}
\newcommand{\cf}{\mathbf{1}}

\begin{document}

\title[Decomposition of the wavelet representation]{A direct
integral decomposition\linebreak[1] of the wavelet representation}
\author[L.-H. Lim]{Lek-Heng Lim}
\address{Department of Mathematics, Malott Hall, Cornell University,
Ithaca, New York 14853-4201}
\email{lekheng@math.cornell.edu}
\curraddr{Department of Pure Mathematics and Mathematical Statistics,
University of Cambridge, Wilberforce Road, Cambridge CB3 0WA, United Kingdom}
\thanks{The third author was supported in part by a grant from NSERC
Canada.}
\author[J. A. Packer]{Judith A. Packer}
\address{Department of Mathematics, National University of Singapore, 10
Kent Ridge Crescent, Singapore 119260}
\email{matjpj@leonis.nus.edu.sg}
\author[K. F. Taylor]{Keith F. Taylor}
\address{Department of Mathematics and Statistics, University of
Saskatchewan, 106 Wiggins Road, Saskatoon, Saskatchewan, Canada S7N 5E6}
\email{taylor@math.usask.ca}
\commby{David R.\ Larson}
\date{November 15, 1999 and, in revised form, February 24, 2000}

\begin{abstract}
In this paper we use the concept of wavelet sets, as introduced by X.\ Dai
and D.\ Larson, to decompose the wavelet representation of the discrete
group associated to an arbitrary $n \times n$ integer dilation matrix as a
direct integral of irreducible monomial representations. In so doing we
generalize a result of F.~Martin and A.~Valette in which they show that
the wavelet representation is weakly equivalent to the regular
representation for the Baumslag-Solitar groups.
\end{abstract}

\maketitle

\section{Introduction}

Several years ago, X.~Dai and D.~Larson introduced an operator algebraic
approach to the study of wavelets, which has proved very useful in the
abstract study of both wavelets and frames \cite{DL1}. In particular,
their study of local commutant and wavelet sets has provided a great
deal of impetus for the research of themselves and others
\cite{DL1,DLS1,DLS2}.

Also, L.~Baggett, both on his own and together with K. Merrill and other
collaborators \cite{B1,BCMO1,BMM1}, has promoted a
representation-theoretic point of view in the study of wavelets and other
types of orthonormal bases of $L^2(\Rn)$ associated to discrete groups. In
\cite{B1}, he decomposed the Stone-von Neumann representation of the
discrete Heisenberg group on $L^2(\R)$ into a direct integral of
representations, and in so doing determined whether or not the translates
of Gabor functions parameterized by certain scales spanned a dense
subspace of $L^2(\R)$.

It is our intention in this paper to use the combined methods of Dai and
Larson \cite{DL1} and Baggett \cite{B1} to further study the notion of
wavelet sets in $\Rn$. Recall from \cite{DLS1} that if $A \in M(n,\Z) \cap
\operatorname{GL}(n,\mathbb{Q})$ satisfies the condition that all of its
eigenvalues have modulus strictly greater than one, it is said to be a
dilation matrix. We then can associate to $A$ a unitary dilation operator
\begin{equation*}
D_A f(t) = \abs{\det A}^{1/2}f(At) \text{,}\displaybreak
\end{equation*}
and to each $v\in \Zn$ a unitary translation operator
\[
T_v f(t) = f(t-v) \text{,}
\]
defined for every $f\in L^2(\Rn)$.

Recall that $\psi \in L^2(\Rn)$ is said to be a \textbf{wavelet} for the
family of unitary operators $\{D_A^m, T_v \mid m\in \Z, v\in \Zn\}$ if
the family of functions
\[
\{D_A^mT_v\psi \mid m\in \Z, v\in \Zn \}
\]
forms an orthonormal basis for $L^2(\Rn)$. Following Dai and Larson, we
define a \textbf{wavelet set} $E$ with respect to the dilation matrix $A$
to be a subset $E \subseteq \Rn$ such that the normalization in $L^2(\Rn)$
of the characteristic function of $E$, $\cf _E$, is the Fourier transform
of a wavelet with respect to the translation and dilation operators
corresponding to $A$. Dai, Larson and Speegle have shown in \cite{DLS1}
that fixing an arbitrary $n \times n$ integer dilation matrix $A$, wavelet
sets with respect to $A$ always exist, and that we necessarily have
$\mu(E)=(2\pi )^n$, for $\mu$ the normalized Lebesgue measure on $\Rn$.

Let $\mathcal{G}(D_A,T_v)$ denote the group in $\mathcal{U}(L^2(\Rn))$
generated by the family of operators $\{D_A, T_v \mid v\in \Zn\}$. Let
$\{e_i \mid 1\leq i\leq n\}$ denote the standard generators for $\Zn$, so
that $(e_i)_j = \delta_{i,j}$. Then the following standard commutation
relations hold:
\[
T_{e_i}D_A = D_AT_{Ae_i} = D_AT_{\sum_{j=1}^{n} a_{j,i}e_j} \text{.}
\]
So $\mathcal{G}(D_A,T_v)$ can be viewed in various ways.  First, it is a
subgroup of $\mathcal{U}(L^2(\Rn))$. Second, it can be viewed as a
unitary representation of the discrete group having $n+1$ generators and
the following $n + n(n-1)/2$ relations:
\[
\left\{
\arraycolsep 0.14em
\begin{array}{rcll}
T_{e_i}D_A &=& D_A\prod_{j=1}^n T_{e_n}^{a_{j,i}}\text{,} &\qquad 1\leq
i\leq n \text{,} \\
T_{e_i}T_{e_j} &=& T_{e_j}T_{e_i}\text{,} &\qquad 1\leq i<j\leq n \text{.}
\end{array}
\right.
\]
We note that groups with generators satisfying the relations above
contain the Baumslag-Solitar groups as a special case; the
$C^{\ast}$-algebras associated to such groups were first observed to be
related to the theory of wavelets by B.~Brenken in Section 6 of
\cite{Br1}. Finally, letting $\Q$ denote the subgroup of $\mathbb{Q}^n$
defined by
\[
\Q = \bigcup_{j=0}^{\infty}\{A^{-j}v \mid v \in \Zn \} \text{,}
\]
$\mathcal{G}(D_A,T_v)$ can be viewed as the image of a representation
of the semidirect product of
$\Q$ and $\Z$, which we will write as $\QZ$. Here for $m\in \Z$,
$\vartheta (m)$ is the automorphism of $\Q$ defined by
\[
\vartheta (m)\beta = A^{-m}\beta
\]
for $\beta \in \Q$, and the group operation in the semidirect product is
defined by
\[
(\beta_1,m_1)\cdot (\beta_2,m_2) = (\beta_1+\vartheta (m_1)\beta_2,
m_1+m_2)
\]
for $(\beta_i,m_i)\in \QZ$. One easily checks that the correspondence
$(\beta ,m) \mapsto T_{\beta}D_A^m$ gives a group isomorphism of $\QZ$
into $\mathcal{G}(D_A,T_v) \subseteq \mathcal{U}(L^2(\Rn))$. Thus we have
a unitary representation of the discrete group $\QZ$ which we term the
\textbf{wavelet representation} of $\QZ$. We use the following notation
for the wavelet representation: for $(\beta,m)\in \QZ$, let
\begin{equation}
\label{eqn W1}
W(\beta,m)=T_{\beta} D_A^m \text{.}
\end{equation}
Then by our remarks above, $W:\QZ \to \mathcal{U}(L^2(\Rn))$ is just the
wavelet representation of $\QZ$.

In this paper, we use the wavelet sets of Dai, Larson and Speegle and
modify the techniques of Baggett discussed above to decompose the wavelet
representation described above as a direct integral of irreducible
representations of $\QZ$. In so doing we are able to extend and generalize
a result of F.~Martin and A.~Valette \cite{MV}, in which they showed that
certain representations of the Baumslag-Solitar groups, of which the
wavelet representation by the univariate wavelet group associated to
dilation by $2$ is a characteristic example, are weakly equivalent to the
regular representations of these groups, hence faithful. Our method will
use the existence of the wavelet sets as proved in \cite{DLS1} to
decompose the wavelet representation as a direct integral of
representations over the wavelet set itself. The main work in the proof
then comes as identifying the representations in this direct integral as
monomial representations, that is, representations induced from characters
on the normal subgroup $\Q$ of $\QZ$.

This paper is an offshoot of the first author's M.Sc.\ thesis done at the
National University of Singapore \cite{LHL}. The third author wishes to
thank the Department of Mathematics and the Wavelets Group at the National
University of Singapore for their hospitality and support during his
visit. Finally, we are grateful to Lawrence Baggett for useful
conversations on the proof of Theorem \ref{thm 3}.
\section{Using the wavelet set to decompose the wavelet representation}

Let $A$ be a fixed $n \times n$ integer dilation matrix.  We first follow
the lead of Dai and Larson \cite{DL1} and study the operators
$\{D_A,T_v\}$ in the Fourier domain. Let $\mathcal{F}$ denote the Fourier
transform on $L^2(\Rn)$, i.e.,
\[
\mathcal{F}f(\xi)=(2\pi)^{-n/2}\int_{\Rn}f(t) e^{-i\ip{t}{\xi}}\, dt
\]
for $f\in L^1\cap L^2(\Rn)$, $t,\xi \in \Rn$, and extending to $L^2(\Rn)$
(which contains $L^1\cap L^2(\Rn)$ as a dense subset) by the usual
limiting process. Hence for $f\in L^2(\Rn)$, $\xi\in \Rn$, we may define
\[
\widehat{D}_A f(\xi) = \mathcal{F}D_A \mathcal{F}^{-1}f(\xi) = \abs{\det
B}^{-1/2}f(B^{-1}\xi )\text{,}
\]
where $B = A^T$, the matrix transpose of $A$, and
\[
\widehat{T}_{\beta} f(\xi) = \mathcal{F}T_{\beta} \mathcal{F}^{-1}f(\xi)
= e^{-i\ip{\beta}{\xi}}f(\xi)\text{,}
\]
where $\beta\in \Q$. As in \cite{DL1}, we calculate the commutant of
$\widehat{D}_A$ and $\{\widehat{T}_v \mid v\in \Zn\}$ in
$\mathcal{B}(L^2(\Rn))$:
\begin{equation}
\label{eqn CVN}
\{\widehat{D}_A, \widehat{T}_v \mid v\in \Zn \}' = \{M_g \mid g\in
L^{\infty}(\Rn), g(\xi)=g(B\xi) \text{ a.e.}, \ \xi\in \Rn\} \text{,}
\end{equation}
where $M_g$ is the multiplication operator defined by
\[
M_g f(\xi)=g(\xi)f(\xi)
\]
for $f\in L^2(\Rn)$, $\xi \in \Rn$. One obtains \eqref{eqn CVN} as
follows: modifying the proof in the one-variable case found in \cite{DL1},
we have
\[
\{\widehat{D}_A,\widehat{T}_v \mid v\in \Zn\}' \subseteq
\{\widehat{T}_{\beta} \mid \beta \in \Q\}' = \{M_g \mid g\in
L^{\infty}(\Rn)\} \text{,}
\]
and computing $\{\widehat{D}_A\}'\cap\{M_g \mid g\in L^{\infty}(\Rn)\}$,
one obtains \eqref{eqn CVN}.

Let $E\subseteq \Rn$ be a wavelet set for the dilation associated to $A$,
so that
\[
\mathcal{F}^{-1}\biggl(\frac{1}{\sqrt{\mu(E)}}\cf_E \biggr)
\]
is a wavelet for $\{D_A, T_v\}$, or to use the definition of Dai, Larson
and Speegle \cite{DLS1}, let $E$ be a \textbf{transformation wavelet set}
for the matrix $B=A^T$, so that the following conditions are satisfied:
\begin{enumerate}
\item  $E$ is a measurable subset of $\Rn$, with $B^j(E) \cap
B^k(E)=\emptyset$, $j \neq k \in \Z$;
\item  $\mu(\Rn \backslash \bigcup_{j\in \Z} B^j(E))=0$;
\item  $E$ is translation congruent to the set $[-\pi,\pi)^n$ modulo the
lattice $\{ 2\pi v \mid v\in \Zn\}$.
\end{enumerate}
If $E$ satisfies conditions (i) and (ii), we say the $E$ \textbf{tiles}
$\Rn$ \textbf{under dilation} by $B$ (cf.\ \cite{DLS2}, where such sets
are said to be \textbf{dilation-congruent} to $[-\pi,\pi)^n$).  The
prototype wavelet set is the subset $E=[-2\pi,-\pi)\cup [\pi,2\pi)$ of
$\R$ corresponding to dilation by $2$ in $\R$ (see \cite{DL1}).

Generalizing the earlier results in \cite{DL1}, Dai, Larson and
Speegle have shown in \cite{DLS1,DLS2} that if $A$ is a dilation matrix,
then $E$ is a wavelet set for the dilation associated to $A$ if and only
if $E$ is a transformation wavelet set for the matrix $B=A^T$, and that
transformation wavelet sets always exist for an arbitrary $n \times n$
dilation matrices with real-valued entries. We concentrate here on the
special case where $A$ has integer entries. It is evident from condition
(iii) for transformation matrix sets that any wavelet set $E$ will satisfy
$\mu(E)=(2\pi)^n$.

We now are able to state our main theorem:
\begin{theorem}
\label{thm 1}
Let $A$ be an $n \times n$ dilation matrix in $M(n,\Z) \cap
\operatorname{GL}(n,\mathbb{Q})$, and let $E\subseteq \Rn$ be a
transformation wavelet set for the matrix $B=A^T$. Then the wavelet
representation $W$ of $\QZ$ defined in \eqref{eqn W1} is unitarily
equivalent to a direct integral of representations $\{\widetilde{W}_x \mid
x\in E \}$, where each $\widetilde{W}_x$ is an irreducible monomial
representation induced from a character on the normal abelian subgroup $\Q
\rtimes_{\vartheta}\{0\}$ of $\QZ$. Indeed, for each $x\in E$, the
representation $\widetilde{W}_x : \QZ \to \mathcal{U}(l^2(\Z))$ is
defined by
\begin{equation}
\label{eqn W2}
[\widetilde{W}_x (\beta,m)g](k)=e^{-i\ip{x}{A^k\beta}}g(k-m)
\end{equation}
for $x\in E$, $k\in \Z$ and $g\in l^2(\Z)$.
\end{theorem}

We will prove Theorem \ref{thm 1} in several stages. We first construct
another Hilbert space isomorphism, which will make the situation more
transparent. Since the sets $B^k(E)$, $k\in \Z$, tile $\Rn$ a.e., we can
define measurable maps
\[
\pi:\Rn\backslash \{0\} \to E \qquad \text{and} \qquad
p:\Rn\backslash \{0\} \to \Z
\]
almost everywhere by the formula
\[
B^{-p(\xi)}\pi(\xi)=\xi \text{,}
\]
i.e., $\pi(\xi)$ is the unique element of $E$ which is $B$-dilation
congruent to $\xi$, and $p(\xi)$ is the unique element of $\Z$ such that
$B^{p(\xi)}\xi \in E$. For example, if $\xi \in E$, $\pi(\xi)=\xi$ and
$p(\xi)=0$, if $\xi\in B^{-1}(E)$, $\pi(\xi)=B\xi$ and $p(\xi)=1$, etc.

Having thus defined the maps $\pi$ and $p$ a.e.\ on $\Rn$, we now define a
measurable transformation $\varphi :\Rn \to E \times \Z$ a.e.\
on $\Rn$ by
\[
\varphi(\xi)=(\pi(\xi),p(\xi)) \text{.}
\]
Finally, we let $\Phi : L^2(\Rn) \to L^2(E \times \Z)$ be the Hilbert
space isomorphism defined by
\[
\Phi f(x,k) = \abs{\det B}^{k/2} f(B^k x)
\]
for $f\in L^2(E\times \Z)$. One readily checks that $\Phi$ preserves
norms and its inverse $\Phi^{-1} :L^2(E \times \Z) \to L^2(\Rn)$ is
defined by
\[
\Phi^{-1}f(\xi)=\abs{\det B}^{p(\xi)/2} f(\pi(\xi),-p(\xi))
\]
for $f \in L^2(\Rn)$ and $\xi \in \Rn$, so that $\Phi$ and $\Phi^{-1}$
are Hilbert space isomorphisms.

We now define maps $\widetilde{D}_A = \Phi \widehat{D}_A \Phi^{-1}$,
and $\widetilde{T}_{\beta } = \Phi \widehat{T}_{\beta } \Phi ^{-1}$,
$\beta \in \Q$, mapping the Hilbert space $L^2(E \times \Z)$ to
itself, and calculate
\begin{align*}
\widetilde{D}_A f(x,k)
&= \Phi (\widehat{D}_A (\Phi^{-1}f))(x,k) \\
&= \abs{\det B}^{k/2} \widehat{D}_A (\Phi^{-1}f)(B^k x) \\
&= \abs{\det B}^{k/2} \abs{\det B}^{-k/2} f(x,k-1) \\
&= f(x,k-1)
\end{align*}
and
\begin{align*}
\widetilde{T}_{\beta} f(x,k)
&= \Phi (\widehat{T}_{\beta} (\Phi^{-1}f))(x,k) \\
&= \abs{\det B}^{k/2}\widehat{T}_{\beta}(\Phi^{-1}f)(B^k x) \\
&= \abs{\det B}^{k/2} e^{-i\ip{\beta}{B^k x}}\abs{\det B}^{-k/2}f(x,k) \\
&= e^{-i\ip{\beta}{B^k x}} f(x,k) \\
&= e^{-i\ip{x}{A^k\beta}} f(x,k) \text{,}
\end{align*}
where $f \in L^2(E \times \Z)$.

We have thus constructed a new unitary representation of the wavelet group
$\QZ$ on the Hilbert space $L^2(E \times \Z)$ which is unitarily
equivalent to the representation $W$ defined in \eqref{eqn W1}. We denote
this representation as follows:
\begin{equation}
\label{eqn W3}
\widetilde{W}(\beta, m)=\widetilde{T}_{\beta} \widetilde{D}_A^m
\end{equation}
for $(\beta,m)\in \QZ$. One computes that
\begin{equation}
\label{eqn W4}
[\widetilde{W}(\beta,m)f](x,k) = e^{-i\ip{x}{A^k\beta}}f(x,k-m)
\end{equation}
for $f\in L^2(E\times \Z)$ and $(\beta,m) \in \QZ$. For fixed $x\in E$,
the formula in \eqref{eqn W4} is given exactly by the formula for the
representation $\widetilde{W}_x$ given in the statement of Theorem
\ref{thm 1}. We now prove that the representations $\widetilde{W}_x$ are
monomial.
\begin{proposition}
\label{prop 1}
Let $E$ be as in the statement of Theorem \ref{thm 1}, and for each $x\in
E$ let $\widetilde{W}_x$ be the representation of $\QZ$ on the Hilbert
space $l^2(\Z)$ defined by
\[
[\widetilde{W}_x(\beta,m)g](k)=e^{-i\ip{x}{A^k\beta}}g(k-m) \text{.}
\]
Then for every $x\in E$, $\widetilde{W}_x$ is an irreducible monomial
representation of $\QZ$, that is, each $\widetilde{W}_x$ is an irreducible
representation which is induced from the one-dimensional representation
$\chi_x$ on the normal abelian subgroup $\Q\rtimes_{\vartheta}\{0\} \cong
\Q$ defined by
\[
\chi_x(\beta)=e^{-i\ip{x}{\beta}}
\]
for $\beta\in \Q$.
\end{proposition}
\begin{proof}
We first note that we have inclusions $\Q \subseteq \mathbb{Q}^n \subseteq
\Rn$ which give us an inclusion monomorphism $\iota: \Q \to \Rn$.  Hence
by Pontryagin duality we have a dual homomorphism $\widehat{\iota}:
\widehat{\R}^n \cong \Rn \to \widehat{\mathbb{Q}}_A$, where the
identification between $\Rn$ and $\widehat{\R}^n$ is given by $\xi \mapsto
\chi_{\xi}(\cdot) = e^{-i\ip{\xi}{\cdot}}$, and we easily check that for
$x\in \Rn$, we have $\widehat{\iota}(x) = \chi_x$ as defined in the
statement of the proposition. Now the action $\vartheta$ of $\Z$ on $\Q$
via automorphisms again via duality gives us a dual action
$\widehat{\vartheta}$ of $\Z$ on $\widehat{\mathbb{Q}}_A$, defined by
$\widehat{\vartheta}(\chi)=\chi \circ \vartheta$, $\chi\in
\widehat{\mathbb{Q}}_A$. By the theory of group representations for
semi-direct product groups (for a reference, see \cite{M1}), it is known
that for every $\chi \in \widehat{\mathbb{Q}}_A$, and for every $m\in \Z$,
the induced representations $\Ind_{\Q}^{\QZ} (\chi)$ and
$\Ind_{\Q}^{\QZ}(\widehat{\vartheta}^m(\chi))$ are equivalent, and that
for fixed $\chi\in \widehat{\mathbb{Q}}_A$, $\Ind_{\Q}^{\QZ} (\chi)$ is
irreducible if and only if $\widehat{\vartheta}^m(\chi) \neq \chi$ for
every $m\in \Z$. Now in our case, for $x\in E$ and $m\in \Z$, we have
\[
[\widehat{\vartheta}^m (\chi_x)](\beta) =\chi_x(\vartheta^m(\beta)) =
e^{-i\ip{x}{A^{-m}\beta}} = e^{-i\ip{B^{-m}x}{\beta}} =
\chi_{B^{-m}x}(\beta) \text{.}
\]
Since $x\in E$, which is a transformation wavelet set for $B$, we know that
the sets $\{B^m(E) \mid m \in \Z \}$ are pairwise disjoint. Hence $B^{-m}x
\neq x$ for every $m\in \Z$. It follows that for every $x\in E$,
the representation $\Ind_{\Q}^{\QZ}(\chi_x)$ is irreducible.

We now show that $\Ind_{\Q}^{\QZ}(\chi_x)$ is equivalent to the
representation $\widetilde{W}_x$ given in the statement of Theorem
\ref{thm 1}. Here we use the construction of induced representations given
in Chapter X of \cite{K1}. In this construction, for a locally compact
group $G$ with closed subgroup $H$, one first needs a Borel cross section
$c:H\backslash G \to G$. Here, $G=\QZ$, $H=\Q \rtimes_{\vartheta}
\{0\}\cong\Q$, and the (right) coset space $H\backslash G \cong\Z$.
Furthermore the groups are discrete, and our cross section $c$ is actually
a group isomorphism defined by
\[
c(k)=(0,k)\in \QZ
\]
for $k \in \Z$. We now use $c$ to define a one-cocycle $\omega$ for the
right action of $G$ on the coset space $H\backslash G$ in the standard
fashion: $\omega: H\backslash G \times G \to H\cong \Q$ is defined by
\begin{align*}
\omega(k,(\beta,m))
&= c(k)(\beta,m)[c(k\cdot (\beta,m))]^{-1} \\
&= (0,k)(\beta,m)(0,-k-m) \\
&= (A^{-k}\beta ,k+m)(0,-k-m) \\
&= (A^{-k}\beta ,0)
\end{align*}
for $k\in H\backslash G \cong \Z$ and $(\beta,m)\in \QZ$. Then we recall
from Chapter X of \cite{K1} that given a unitary representation $\Lambda$
of the group $H$ on the Hilbert space $\mathcal{H}$, the induced
representation $\Ind_H^G(\Lambda)$ has as its representation space
$L^2(H\backslash G)\otimes \mathcal{H}\cong L^2(H\backslash G,
\mathcal{H})$ where the measure on $H\backslash G$ is quasi-invariant
under translation by $G$. If this measure is invariant, then the formula
for $\Ind_H^G(\Lambda)$ is given by
\[
[\Ind_H^G(\Lambda)(x)f](\bar{y}) =[\Lambda(\omega(\bar{y},x))f](\bar{y}x)
\]
for $x\in G$, $\bar{y}\in H\backslash G$, $f\in L^2(H\backslash G,
\mathcal{H})$. Here all our groups are discrete, and the
translation-invariant measure on the coset space $H\backslash G\cong \Z$
is just the counting measure.  Hence for $x\in E$ and $\chi_x$ defined as
in the statement of Proposition \ref{prop 1} we have
$\mathcal{H}=\mathbb{C}$ and $L^2(H\backslash G,\mathcal{H})\cong
l^2(\Z),$ so that with respect to these identifications,
\begin{align*}
[\Ind_{\Q}^{\QZ}(\chi_x)(\beta,m)f](k)
&= \chi_x(\omega(k,(\beta,m)))f(k\cdot (\beta,m)) \\
&= e^{-i\ip{x}{A^{-k}\beta}}f(k+m)
\end{align*}
for $k \in \Z$, $(\beta,m)\in \QZ$ and $f\in l^2(\Z)$. Now this is not
exactly the same as our formula for $\widetilde{W}_x$, but defining the
unitary involution $V:l^2(\Z)\to l^2(\Z)$ by the formula $Vf(k)=f(-k)$,
one easily checks that
\begin{align*}
V^{-1}[\widetilde{W}_x(\beta,m)]V &= \Ind_{\Q}^{\QZ}(\chi_x)(\beta,m)
\end{align*}
for $(\beta,m)\in\QZ$ and $x\in E$. This completes the proof of the
proposition.
\end{proof}

We are now ready to finish the proof of Theorem \ref{thm 1}.
\begin{proof}[Proof of Theorem \ref{thm 1}]
We have constructed an equivalence between the wavelet representation $W$
of \eqref{eqn W1} and the representation $\widetilde{W}$ of \eqref{eqn
W3}, so it only remains to show that $\widetilde{W}$ is a direct integral
of irreducible monomial representations. Now in general the theory of
direct integrals of measurable fields of Hilbert spaces and direct integrals
of unitary representations is very technical
(see, for example, \cite{M1} for a reference). But in our situation, the
representation space for the representation $\widetilde{W}$ is
$L^2(E\times \Z)$ which can also be viewed as $L^2(E)\otimes l^2(\Z)$, and
there are no technicalities involved in showing that this last Hilbert
space is exactly the direct integral $\int_E^{\oplus} (l^2(\Z))_x\, dx$.
With respect to this decomposition of $L^2(E \times \Z)$ it is clear that
the decomposition of $\widetilde{W}$ over the measurable subset $E$ is
given exactly by the representations $\{ \widetilde{W}_x \mid x\in E \}$
defined in \eqref{eqn W2}. Finally, Proposition \ref{prop 1} has shown
that for each $x\in E$ the representation $\widetilde{W}_x$ is an
irreducible monomial representation of $\QZ$, and the proof of Theorem
\ref{thm 1} is complete.
\end{proof}

\begin{remark}
\label{rem set}
We note here that to prove Theorem \ref{thm 1} we did not need to use all
of the properties (i)--(iii) of transformation wavelet sets. Indeed, all
we really used were properties (i) and (ii): we did not need that $E$ be
$2\pi$-translation congruent to $[-\pi,\pi)^n$. The pairwise disjointness
of the sets $\{ B^j(E) \mid j\in \Z \}$ and the fact that $\mu(\Rn
\backslash \bigcup_{j\in \Z}~B^j(E))=0$ allowed us to construct the
isomorphism of measure spaces $\varphi: \Rn \to E\times \Z$ which was
crucial in the direct integral procedure. Property (ii) ensures that the
representations $\widetilde{W}_x$ are irreducible.  It is only when one
wants to show that $\mathcal{F}^{-1}(\cf_E/\sqrt{\mu(E)})$ is a wavelet
for the dilation system $\{D_A, T_v \mid v\in \Zn \}$ that property (iii)
becomes essential. The next corollary gives a proof of this result,
originally due to Dai and Larson (\cite{DL1}, cf.\ Lemmas 2.2 and 2.3)
which uses the framework that we have set up. The main difference between
our method and that of Dai and Larson is the slightly more transparent
presentation given to the operators $\{ \widetilde{D}_A^m \mid m\in \Z \}$
as operators which are just translations in the second variable on our
realization of the Hilbert space as $L^2(E\times \Z)$.
\end{remark}

\begin{corollary}
\label{cor 1}
Let $A$ be a dilation matrix in $M(n,\Z) \cap \operatorname{GL}
(n,\mathbb{Q})$ with $A^T=B$, and suppose that $E$ is a subset of $\Rn$
with finite measure such that the sets $\{B^j(E) \mid j \in \Z \}$ are
pairwise disjoint and $\mu(\Rn \backslash \bigcup_{j\in \Z}B^j(E))=0$.
Then $\psi = \mathcal{F}^{-1}(\cf_E/\sqrt{\mu(E)})$ is a wavelet for the
system $\{D_A, T_v \mid v\in \Zn \}$ if and only if the sets $\{ E+2\pi v
\mid v\in \Zn \}$ are pairwise disjoint up to sets of measure zero and
$\mu(\Rn \backslash \bigcup_{v\in \Zn}(E+2\pi v))=0$.
\end{corollary}

\begin{proof}
The content of Remark \ref{rem set} shows that under the hypotheses of the
Corollary, the proof of Theorem \ref{thm 1} still works. Keeping the
notation of Theorem \ref{thm 1} and Proposition \ref{prop 1}, we see that
the set $\{ D_A^m T_v\psi \mid m\in \Z, v\in \Zn\}$ will be an orthonormal
basis for $L^2(\Rn)$ if and only if the set $\{\widetilde{D}_A^m
\widetilde{T}_v \Phi^{-1}(\cf_E/\sqrt{\mu(E)}) \mid m\in \Z, v\in \Zn \}$
is an orthonormal basis for $L^2(E\times \Z)$. Now
\[
\Phi^{-1}\biggl(\frac{1}{\sqrt{\mu(E)}}\cf_E \biggr)
= \frac{1}{\sqrt{\mu(E)}} \cf_{E\times\{0\}} \text{,}
\]
where $E\times \{0\}$ is viewed as a subset of $E\times \Z$, and since
$\widetilde{D}_A^m f(x,k)=f(x,k-m)$ it is clear that for any $f\in
L^2(E\times \Z)$ which is supported on $E\times \{0\}$, the function
$\widetilde{D}_A^m f$ will be supported on $E\times \{m\}$, $m\in \Z$. It
follows from this and an inspection of the formula for $\widetilde{D}_A$
that if $\{f_i \mid i\in \text{index set }I \}$ is a set of orthonormal
functions whose support lies in $E\times \{0\}$ which form a basis for
$L^2(E\times \{0\})$, then the set $\{\widetilde{D}_A^m f_i \mid m\in \Z,
i \in I \}$ will form an orthonormal basis for $L^2(E\times \Z)$. Thus in
order to obtain the results of the corollary it is enough to show that the
stated conditions on $E$ are necessary and sufficient to ensure that the
set $\{ \widetilde{T}_v (\cf_{E\times \{0\}}/\sqrt{\mu(E)}) \mid v\in
\Zn \}$ is an orthonormal basis for $L^2(E\times \{0\})$. Recall that
for $v \in \Zn$ we have
\[
\widetilde{T}_v \biggl(\frac{1}{\sqrt{\mu(E)}} \cf_{E\times \{0\}} \biggr)
(x,0) = \frac{1}{\sqrt{\mu(E)}} e^{-i\ip{v}{x}}
\]
for $x\in E$ and it is well known (cf.\ \cite{DL1} Lemma 4.2,
Lemma 4.3, \cite{DLS2}) that as $v$
varies over all of $\Zn$, the functions above will give an orthonormal
basis for $L^2(E)$ if and only if $E$ is translation congruent to
$[-\pi,\pi)^n$, which is equivalent to the conditions on $E$ given in the
statement of the corollary.
\end{proof}
\section{Weak equivalence of the wavelet representation\protect\linebreak[1] and the regular
representation}

Let $d\in \mathbb{N}$, and for $d \geq 2$ let $BS_d$ denote the
Baumslag-Solitar groups considered by F.~Martin and A.~Valette \cite{MV};
we recall that $BS_d$ has two generators $a$ and $b$ satisfying the
relation $aba^{-1}=b^d$. It is clear that with respect to our notation,
$BS_d$ is the wavelet group corresponding to the dilation $A=(d)$ on
$\R^1$. In Theorem 8 of \cite{MV}, it is shown that the regular
representation of $BS_d$ on $l^2(BS_d)$ and the wavelet representation of
$BS_d$ on $L^2(\R)$ are weakly equivalent. The proof of our Theorem
\ref{thm 1} allows us to extend this result to our wider class of
semidirect product wavelet groups $\QZ$, where as usual for $n\geq 1$, $A$
is a dilation matrix in $M(n,\Z) \cap \operatorname{GL}(n,\mathbb{Q})$,
$\Q$ is the subgroup $\{A^m v \mid m\in \Z, v\in \Zn\}$ of $\mathbb{Q}^n$
determined by $A$, and $\vartheta$ is the automorphism of $\Q$ associated
to $A$. We thank Lawrence Baggett for suggesting that we use the reference
\cite{F2}, which allowed us to considerably simplify our original proof.
\begin{theorem}
\label{thm 3}
Let $n\in \mathbb{N}$ and let $A$ be a dilation matrix in $M(n,\Z) \cap
\operatorname{GL}(n,\mathbb{Q})$. Let $\QZ$ be the wavelet group
associated to these parameters described in Section $2$. Then the regular
representation and the wavelet representation of $\QZ$ are weakly
equivalent. Hence the wavelet representation of $\QZ$ is faithful.
\end{theorem}

\begin{proof}
 Let $E$ be a transformation wavelet set for $B=A^T$. The proof of
Theorem \ref{thm 1} shows that the wavelet representation $W$ is unitarily
equivalent to a representation $\widetilde{W}$ which can be expressed as a
direct integral
\[
\int_E^{\oplus} \widetilde{W}_x\, dx \text{,}
\]
where for each $x\in E$, $\widetilde{W}_x$ is an irreducible representation
which is unitarily equivalent to $\Ind_{\Q}^{\QZ}(\chi_x)$ defined on the
Hilbert space $l^2(\Z)$. We also noted in the proof of Proposition
\ref{prop 1} that for $x\in E$ and $m\in \Z$ the representations
$\Ind_{\Q}^{\QZ}(\chi_x)$ and $\Ind_{\Q}^{\QZ}(\chi_{B^m x})$ are
equivalent to one another. Hence the representation $\widetilde{W}$, and
consequently the wavelet representation $W$, is weakly contained in the set of all of the
monomial representations
\[
\bigl\{\Ind_{\Q}^{\QZ}(\chi_y) \bigm| y\in \bigcup_{j\in \Z} B^j(E)
\bigr\} \text{.}
\]
Let us now consider the regular representation $R$ of $\QZ$. Recall that
$R$ itself is an induced representation,
\[
R=\Ind_{\{1\}}^{\QZ}(1) \text{,}
\]
where here $\{ 1\}$ represents the trivial subgroup of $\QZ$, and the
representation of $\{ 1\}$ we are inducing is the trivial one-dimensional
representation on $\mathbb{C}$. By the theory of induction in stages
\cite{M1}, we can write
\[
R\cong \Ind_{\Q}^{\QZ}\left[ \Ind_{\{1\}}^{\Q}(1)\right] \text{.}
\]
Now let $R_1$ denote the regular representation of the discrete group $\Q$,
that is,
\[
R_1=\Ind_{\{1\}}^{\Q}(1) \text{.}
\]
By the Fourier theory connecting compact and discrete abelian groups,
$R_1$ is equivalent to a representation on the space
$L^2(\widehat{\mathbb{Q}}_A)$ and in fact can be represented as a direct
integral of characters as follows.
Setting ${\widehat{\mathbb Q}_A}=\Sigma_A,$ we have
\[
R_1\cong \int_{\Sigma_A}^{\oplus}\gamma \, d\gamma \text{,}
\]
where each $\gamma \in \Sigma_A$, being a character of ${\mathbb{Q}}_A,$ is
exactly a one-dimensional representation on the space $\mathbb{C}$. We
thus obtain that $R$ is equivalent to
\[
\Ind_{\Q}^{\QZ}\left[\displaystyle{\int_{\Sigma_A}^{\oplus}\gamma
\, d\gamma }\right] \text{,}
\]
and since the processes of taking direct integrals of representations and inducing from a subgroup commute,
we have that
\[
R\cong
\int_{\Sigma_A}^{\oplus}\left[\Ind_{\Q}^{\QZ}(\gamma
)\right]\, d\gamma,
\]
the representation on the right-hand side being defined on the Hilbert space
\[
\int_{\Sigma_A}^{\oplus}\left[ l^2(\Z)\right]_{\gamma }\,
d\gamma \text{.}
\]
Since we have already shown in the proof of Proposition \ref{prop 1} that
\[
\bigl\{\chi_y \bigm| y\in \textstyle{\bigcup_{j \in \Z }}B^j(E)\bigr\}
\subseteq \bigl\{\chi_t \bigm| t\in \Rn \bigr\} \subseteq \Sigma_A =
\widehat{\mathbb{Q}}_A \text{,}
\]
we have shown directly that the representation
$\widetilde{W}$ and hence the wavelet representation $W$ is weakly contained in
$R.$ Of course, as
with the Baumslag-Solitar groups studied in \cite{MV}, this also follows
immediately from the fact that our group $\QZ$, being the semi-direct
product of two abelian groups, is amenable.

To show that $R$ is weakly contained in $W$, we first
note that our hypotheses on $B$
together with Proposition 2.4 of \cite{Br1} show that if we let
$\widehat{\iota}: \Rn \to {\Sigma}_A$ be the monomorphism constructed in
the proof of Proposition \ref{prop 1}, then the range of $\widehat{\iota}$
is dense in $\Sigma_A$. Since $\mu(\Rn \backslash \bigcup_{j\in \Z}
B^j(E))=0$, $\bigcup_{j\in \Z}B^j(E)$ is dense in $\Rn$ in the usual
topology, so that $\widehat{\iota}(\bigcup_{j\in \Z} B^j(E))$ is dense in
$\widehat{\iota}(\Rn)$ in the relative topology. It thus follows that
$\widehat{\iota}(\bigcup_{j\in \Z} B^j(E))$ is dense in $\Sigma_A$. Hence
\[
\overline{\widehat{\iota}\bigl(\textstyle{\bigcup_{j\in \Z}B^j(E)} \bigr)}
= \overline{\bigl\{\chi_y \bigm| y\in \textstyle{\bigcup_{j\in \Z}
B^j(E)}\bigr\}} =\Sigma_A \text{.}
\]
The argument that now follows is similar to the proof of Lemma 4 in
\cite{MV}, which was not used in the proof of Theorem 8 of \cite{MV}.
We now use the continuity of the induction process (cf.\ \cite{F2}) to
deduce that
\begin{equation}
\label{eq wk}
\overline{\bigl\{ \Ind_{\Q}^{\QZ}(\chi_y) \bigm| y\in \bigcup_{j\in \Z}
B^j(E)\bigr\}}
= \overline{\bigl\{\Ind_{\Q}^{\QZ}(\gamma) \bigm| \gamma\in \Sigma_A
\bigr\}} \text{,}
\end{equation}
where the closure is taken in the hull-kernel topology. We have already
shown that $W$ is unitarily equivalent to the direct integral of the
representations
\[
\bigl\{ \Ind_{\Q}^{\QZ}(\chi_y) \bigm| y\in E\bigr\} \text{,}
\]
and that for fixed $y\in E$ and $m\in \Z$, $\Ind_{\Q}^{\QZ}(\chi_{B^m y})$
is unitarily equivalent to the representation $\Ind_{\Q}^{\QZ}(\chi_y)$.
Since $R$ is the direct integral of the representations $\Ind_{\Q}^{\QZ}
(\gamma)$, it follows from \eqref{eq wk} that $R$ is weakly contained in $W,$
so that $W$ and $R$ are weakly equivalent.  Finally, as noted in \cite{MV},
$R$ is faithful since $\QZ$ is amenable, and it follows from the definition of
weak equivalence that $W$ is faithful,
thus generalizing Theorem 8 proved in \cite{MV} to the wider class of groups
$\QZ.$
\end{proof}

\end{document}